\documentclass[a4paper, 12pt]{article}
\usepackage{}

\usepackage{mathrsfs}
\usepackage{epsfig}
\usepackage{amsmath}
\usepackage{amssymb}
\usepackage{latexsym}
\usepackage{amsfonts}
\usepackage{amsthm}

\usepackage[numbers,sort&compress]{natbib}
\usepackage{graphicx}
\ifx\pdfoutput\undefined  \else
 \fi

\usepackage[centerlast]{caption2}

\usepackage{color}

\usepackage{txfonts}
\usepackage{amssymb}
\usepackage{mathrsfs}
\usepackage{amsfonts}

\newtheorem{theorem}{Theorem}[section]

\usepackage{latexsym,amssymb}
\usepackage{amsmath}

\newcommand{\F}{{\mathbb F}}

\begin{document}

\title{{Lower bound for the Erd\H{o}s-Burgess constant of finite commutative rings}}
\author{
Guoqing Wang\\
\small{Department of Mathematics, Tianjin Polytechnic University, Tianjin, 300387, P. R. China}\\
\small{Email: gqwang1979@aliyun.com}
\\
}
\date{}
\maketitle

\begin{abstract} Let $R$ be a finite commutative unitary ring. An idempotent in $R$ is an element $e\in R$ with $e^2=e$.  The Erd\H{o}s-Burgess constant associated with the ring $R$ is the smallest positive integer $\ell$ such that for any given $\ell$ elements (repetitions are allowed) of $R$, say $a_1,\ldots,a_{\ell}\in R$, there must exist a nonempty subset $J\subset \{1,2,\ldots,\ell\}$ with $\prod\limits_{j\in J} a_j$ being an idempotent.
In this paper, we give a lower bound of the Erd\H{o}s-Burgess constant in a finite commutative unitary ring in terms of all its maximal ideals, and prove that the lower bound is attained in some cases. The result unifies some recently obtained theorems on this invariant.
\end{abstract}

\noindent{\small {\bf Key Words}: {\sl  Erd\H{o}s-Burgess constant;  Davenport constant;  Zero-sum; Idempotent-product free sequences; Finite commutative rings}}

\section {Introduction}

Let $\mathcal{S}$ be a nonempty commutative semigroup, endowed with a binary associative operation $*$. Let ${\rm E}(\mathcal{S})$ be the set consisting of all idempotents of $\mathcal{S}$, where $e\in \mathcal{S}$ is said to be an idempotent if $e*e=e$.
Idempotent is one of central notions in Semigroup Theory and Algebra.
One of our interest to combinatorial properties concerning idempotents in semigroups comes from a question of P. Erd\H{o}s to D.A. Burgess (see \cite{Burgess69} and \cite{Gillam72}), which is stated as follows.

{\sl Let $\mathcal{S}$ be a finite nonempty semigroup of order $n$. A sequence of terms from $\mathcal{S}$ of length $n$ must contain one or more terms whose product, in some order, is an idempotent?}

Burgess \cite{Burgess69} in 1969 gave an answer to this question in the case when $\mathcal{S}$ is commutative or contains only one idempotent. Note that every nonempty finite semigroup contains at least one idempotent (see \cite{Grillet monograph} Chapter I, Corollary 5.9).
D.W.H. Gillam, T.E. Hall and N.H. Williams \cite{Gillam72} proved that a sequence $T$ over any finite semigroup $\mathcal{S}$ of length at least $|\mathcal{S}\setminus{\rm E}(\mathcal{S})|+1$ must contain one or more terms whose product, in the order induced from the sequence $T$, is an idempotent, and therefore, completely answered Erd\H{o}s' question.
The Gillam-Hall-Williams Theorem was extended to infinite semigroups by the author \cite{wangidempotent} in 2019. It was also remarked  that the bound $|\mathcal{S}\setminus {\rm E}(\mathcal{S})|+1$,  although is
optimal for general semigroups $\mathcal{S}$, can be improved, at least in principle, for specific classed of semigroups.
Naturally, one combinatorial invariant was aroused by Erd\H{o}s' question with respect to idempotents of semigroups. Since we deal with the multiplicative semigroup of a commutative ring in this paper, we introduce only the definition of this invariant for commutative semigroups here.

\noindent \textbf{Definition A.} (\cite{wangidempotent}, Definition 4.1) \ {\sl For a commutative semigroup $\mathcal{S}$, define the {\bf Erd\H{o}s-Burgess constant} of $\mathcal{S}$, denoted by $\textsc{I}(\mathcal{S})$, to be the least $\ell\in \mathbb{N}\cup \{\infty\}$ such that every sequence $T$ of terms from $\mathcal{S}$ and of length $\ell$ must contain one or more terms whose product is an idempotent.}

Note that if the commutative semigroup $\mathcal{S}$ is finite, Gillam-Hall-Williams Theorem definitely tells us that the Erd\H{o}s-Burgess constant of  $\mathcal{S}$ exists, i.e., $\textsc{I}(\mathcal{S})\in \mathbb{N}$ is finite and bounded above by $|\mathcal{S}\setminus{\rm E}(\mathcal{S})|+1$.  In particular,
when the semigroup $\mathcal{S}$ happens to be a finite abelian group, the Erd\H{o}s-Burgess constant reduces to a classical combinatorial invariant, the Davenport constant. The Davenport constant of a finite abelian group $G$, denoted by ${\rm D}(G)$, is defined as the smallest positive integer $\ell$  such that every sequence of terms from $G$ of length $\ell$ must contain one or more terms with the product being the identity element of $G$. This invariant was popularized by H. Davenport in the 1960's, notably for its link with algebraic number theory (as reported in \cite{Olson1}).
For the  progress about ${\rm D}(G)$, the reader may
consult e.g., \cite{GaoGeroldingersurvey, GLP12} on commutative groups, \cite{GaoLiPeng,GDavid} on noncommutative groups, and \cite{Deng, wangDavenportII, wangAddtiveirreducible, wang-zhang-qu, wang-zhang-wang-qu} on commutative semigroups and commutative rings.

The author \cite{wangErdosburgess} obtained some conditions such that the Erd\H{o}s-Burgess constant exists in general commutative rings, which is stated in Theorem B below.

\noindent $\bullet$ {\small For any commutative unitary ring $R$, let ${\rm U}(R)$ be the group of units and $\mathcal{S}_R$ the multiplicative semigroup of the ring $R$.}

\noindent \textbf{Theorem B.} \cite{wangErdosburgess}  \ {\sl Let $R$ be a commutative unitary ring. If $\textsc{I}(\mathcal{S}_R)$ is finite, then one of the following two conditions holds:

\noindent (i) The ring $R$ is finite;

\noindent (ii) The Jacobson radical $J(R)$ is finite and $R\diagup J(R)\cong B \times \prod\limits_{i=1}^t \mathbb{F}_{q_i}$, where  $B$ is an infinite Boolean unitary ring, and $\mathbb{F}_{q_1},\ldots, \mathbb{F}_{q_t}$ are finite fields with $0\leq t\leq \textsc{I}(\mathcal{S}_R)-1$ and prime powers $q_1,\ldots, q_t>2$.}

The above Theorem B asserts that the Erd\H{o}s-Burgess constant exists only for {\bf finite} commutative rings except for a family of infinite commutative rings with very special forms given as (ii) above. That is, to study this invariant in the realm of commutative rings, we may consider it only for {\sl finite} commutative rings.
Recently, J. Hao, H. Wang and L. Zhang obtained a sharp lower bound of the Erd\H{o}s-Burgess constant in two classes of finite commutative rings, i.e.,  in the residue class ring $\mathbb{Z}\diagup n \mathbb{Z}$ and in the quotient ring $\mathbb{F}_q[x]\diagup K$ of the polynomial ring $\mathbb{F}_q[x]$ modulo a nonzero proper ideal $K$, which are stated as Theorem C and Theorem D below.

\noindent \textbf{Theorem C.} \cite{wang-Hao-zhang}  \ {\sl Let $q$ be a prime power, and let $\F_q[x]$ be the ring of polynomials over the finite field $\F_q$. Let $R=\F_q[x]\diagup K$ be a quotient ring of $\F_q[x]$ modulo any nonzero proper ideal $K$. Then ${\rm I}(\mathcal{S}_R)\geq {\rm D}({\rm U}(R))+\Omega(K)-\omega(K),$ where $\Omega(K)$ is the number of the prime ideals (repetitions are counted) and $\omega(K)$ the number of distinct prime ideals
in the factorization when $K$
is factored into a product of prime ideals. Moreover, equality holds for the case when
$K$ is factored into either a power of some prime ideal or a product of some pairwise distinct prime ideals.}

\noindent \textbf{Theorem D.} \cite{HaoWangZhang} \ {\sl Let $n>1$ be an integer, and let $R=\mathbb{Z}\diagup n \mathbb{Z}$ be the ring of integers modulo $n$.
Then
${\rm I}(\mathcal{S}_R)\geq {\rm D}({\rm U}(R))+\Omega(n)-\omega(n),$ where $\Omega(n)$ is the number of primes occurring in the prime-power decomposition of $n$ counted with
multiplicity, and $\omega(n)$ is the number of distinct primes.  Moreover, equality holds if $n$ is a prime power
or a product of pairwise distinct primes.}

Note that both the ring $\F_q[x]$ of polynomials over a finite field and the ring $\mathbb{Z}$ of integers are principal ideal domains and definitely are Dedekind domains, in which every nonzero proper ideal has a unique factorization as a product of prime ideals (which are also maximal ideals).
N. Kravitz and A. Sah \cite{KravitzSah} completely generalized Theorems C and D to the finite quotient ring of any Dedekind domain.  For a Dedekind domain $D$ and a nonzero proper ideal $K$ of $D$, let $\Omega(K)$ be the total number of prime ideals in the prime ideal factorization of $K$ (with multiplicity), and let $\omega(K)$ be the number of distinct prime ideals in this factorization.

\noindent \textbf{Theorem E.} \cite{KravitzSah}  \ {\sl Let $D$ be a Dedekind domain and $K$ a nonzero proper ideal of $D$ such that $R=D\diagup K$ is a finite ring.  Then ${\rm I}(\mathcal{S}_R)\geq {\rm D}({\rm U}(R))+\Omega(K)-\omega(K).$
Moreover, equality holds if $K$ is either a power of a prime ideal or a product of distinct prime ideals.}

In this paper, we shall obtain a sharp lower bound for the Erd\H{o}s-Burgess constant of a general finite commutative unitary ring, which generalizes Theorem E, and definitely deduces
Theorem C and Theorem D as consequences. To give the theorem, we need one notion as follows.

Let $R$ be a finite commutative unitary ring, and let $N$ be an ideal of $R$.  For any nonnegative integer $i$, let $N^i$ be the $i$-th power of the ideal $N$. In particular, we define $N^{0}=R.$ Define the {\bf index} of the ideal $N$, denoted by ${\rm Ind}(N)$, to be the least nonnegative integer $k$ such that $N^k=N^{k+1},$ equivalently, the descending chain of ideals $N^{0}\supsetneq N^{1}\supsetneq  \cdots \supsetneq N^k=N^{k+1}=N^{k+2}=\cdots$ becomes stationary starting from $N^k$. Now we are in a position to give the theorem of this paper.

\begin{theorem}\label{Theorem main}
Let $R$ be a finite commutative unitary ring. Then
$${\rm I}(\mathcal{S}_R)\geq {\rm D}({\rm U}(R))+ \sum\limits_{M} ({\rm Ind}(M)-1)$$ where $M$ is taken over all distinct maximal ideals of $R$.  Moreover, equality holds if $R$ is a local ring or all its maximal ideals have the indices one.
 \end{theorem}

 \noindent {\bf Remark.}  Since any ideal $K$ is prime if and only if it is maximal in a finite commutative unitary ring,  we can restate Theorem \ref{Theorem main} in terms of prime ideals of $R$ as the following theorem. Moreover,
 we will show that Theorem \ref{Theorem main} implies Theorem E, and therefore Theorem C and Theorem D as consequences in the final Concluding remarks section.

\noindent {\bf Theorem.}  Let $R$ be a finite commutative unitary ring. Then
${\rm I}(\mathcal{S}_R)\geq {\rm D}({\rm U}(R))+ \sum\limits_{M} ({\rm Ind}(M)-1)$ where $M$ is taken over all distinct prime ideals of $R$.  Moreover, equality holds if $R$ has a unique prime ideal or all its prime ideals have the indices one.

\section{The Proof}

To prove Theorem \ref{Theorem main}, we need to introduce some necessary notation and terminology, and follow the notation of A. Geroldinger, D.J. Grynkiewicz and others used for sequences over
groups (see \cite{GDavid} for example).
For integers $a,b\in \mathbb{Z}$, we set $[a,b]=\{x\in \mathbb{Z}: a\leq x\leq b\}$.
Let $R$ be a finite commutative unitary ring with multiplication *, an element $e$ in $R$ is called an idempotent if $e*e=e$. A sequence $T$ of terms from $R$ is a multi-set for which the repetitions of terms are allowed, denoted as
\begin{equation}\label{equation sequence T}
T=a_1a_2\cdot\ldots\cdot a_{\ell}=\mathop{\bullet}\limits_{i\in [1,\ell]} a_i
\end{equation}
 where $a_i\in R$ for each $i\in [1,\ell]$. By $\cdot$ we denote the concatenation of sequences which differs from the notation * used for the multiplication of the ring $R$.  By $|T|$ we denote the length of the sequence $T$. In particular, $|T|=\ell$ for the sequence $T$ in \eqref{equation sequence T}.
 Since the ring $R$ is commutative and the order of terms of $T$ does not matter for the combinatorial property of sequences which will be investigated here, then for any permutation $\tau$ of $[1,\ell]$ we always identify $a_{\tau(1)}a_{\tau(2)}\cdot\ldots\cdot a_{\tau(\ell)}$ with the sequence $T$ given in \eqref{equation sequence T}.  For any subset $X\subseteq [1,\ell]$, we say $T'=\mathop{\bullet}\limits_{i\in X} a_i$ is a subsequence of $T$, denoted by $T'\mid T$, and in particular we call $T'$ a {\sl proper} [resp. {\sl nonempty}] subsequence of $T$ if $X\neq [1,\ell]$ [resp. if $X\neq \phi$].  By $\varepsilon$ we denote the
{\sl empty sequence}  with $|\varepsilon|=0$, which is a proper subsequence of every nonempty sequence.
Let $\pi(T)=\prod\limits_{i\in [1,\ell]} a_i=a_1*a_2*\cdots* a_{\ell}$ be the product of all terms of $T$. We adopt the convention that $\pi(\varepsilon)=1_R$.
By $\prod(T)$ we denote the set of elements of $R$ that can be represented as a product of one or more terms from $T$, i.e.,
$$\prod(T)=\{\pi(T'): T' \mbox{ is taken over all nonempty subsequences of }T\}.$$ We call $T$ an {\sl idempotent-product free sequence} provided that no idempotent of $R$ can be represented as a product of one or more terms from $T$, i.e., $\prod(T)$ contains no idempotent.

\medskip

\noindent {\bf Proof  of Theorem \ref{Theorem main}.} \
Let $M_1,\ldots,M_r$ be all maximal ideals of $R$ with indices $k_1, \ldots, k_r$ respectively, where $r\geq 1$. For each $i\in [1,r]$, since $M_i^ {k_i-1}\supsetneq  M_i^ {k_i}$, we can take some element $x_i$ of $M_i^ {k_i-1}\setminus M_i^ {k_i}$. Since $x_i$ is a finite sum of products of the form $a_{1}*a_{2}*\cdots *a_{k_i-1}$ where $a_{1}, a_{2},\ldots, a_{k_i-1}\in M_i$, it follows that there exist $k_i-1$ elements, say $y_{i,1},\ldots,y_{i,{k_i-1}}\in M_i$ such that $$y_{i,1}*\cdots *y_{i,{k_i-1}}\in M_i^{k_i-1}\setminus  M_i^{k_i}.$$
Note that
\begin{equation}\label{equation yij in Mi|Ji|}
\prod\limits_{j\in J_i} y_{i,j}\in M_i^{|J_i|}\setminus M_i^{|J_i|+1} \mbox{ for any } i\in [1,r] \mbox{ and any subset } J_i\subseteq [1,k_i-1].
\end{equation}
Since $M_1^ {k_1}, \ldots, M_r^ {k_r}$ are pairwise coprime ideals of $R$, by the Chinese Remainder Theorem, for any $i\in [1,r]$ and $j\in [1,k_i-1]$,
we can find an element $\widetilde{y}_{i, j}$ of $R$ such that
\begin{equation}\label{equation tildeyij=yij}
\widetilde{y}_{i, j}\equiv y_{i, j} \pmod {M_i^{k_i}}
\end{equation}
and \begin{equation}\label{equation tildeyij=1}
\widetilde{y}_{i, j}\equiv 1_R \pmod {M_t^{k_t}} \mbox{ where } t\in [1,r]\setminus \{i\}.
\end{equation}
By the definition of the Davenport constant, we can take a sequence $V$ of terms from the group ${\rm U}(R)$ with length
\begin{equation}\label{equation |V|_D-1}
|V|={\rm D}({\rm U}(R))-1
\end{equation}
 such that the identity element $1_R$ of the group ${\rm U}(R)$ can not be represented as a product of one or more terms from $V$, i.e.,
\begin{equation}\label{equation 1 notin prod(V)}
1_R\notin \prod(V).
\end{equation}
 Let
 \begin{equation}\label{equation form of the sequence T}
 T=V\cdot (\mathop{\bullet}\limits_{j\in [1,k_1-1]} \widetilde{y}_{1, j})\cdot \ldots\cdot (\mathop{\bullet}\limits_{j\in [1,k_r-1]} \widetilde{y}_{r, j}).
 \end{equation}

\noindent {\bf Assertion.} \  The sequence $T$ is idempotent-product free.

\noindent {\sl Proof.} \
Assume to the contrary that $T$ contains a nonempty subsequence $T'$ such that $\pi(T')$ is an idempotent.
Since $1_R$ is the unique idempotent in ${\rm U}(R)$, it follows from \eqref{equation 1 notin prod(V)} that at least one term of $(\mathop{\bullet}\limits_{j\in [1,k_1-1]} \widetilde{y}_{1, j})\cdot \ldots\cdot (\mathop{\bullet}\limits_{j\in [1,k_r-1]} \widetilde{y}_{r, j})$ appears in the sequence $T'$, then by rearranging the indices $i\in [1,r]$ and $j\in [1, k_i-1]$ we may assume without loss of generality that
\begin{equation}\label{equation T'=}
T'=V'\cdot W' \cdot (\mathop{\bullet}\limits_{j\in [1,n]} \widetilde{y}_{1, j})
\end{equation}
 where $V'\mid V$,
\begin{equation}\label{equation W' is a subsequence of}
W'\mid (\mathop{\bullet}\limits_{j\in [1,k_2-1]} \widetilde{y}_{2, j})\cdot \ldots\cdot (\mathop{\bullet}\limits_{j\in [1,k_r-1]} \widetilde{y}_{r, j})
\end{equation}
 and
 \begin{equation}\label{equation 1 leq t leq k_1-1}
 1\leq n\leq k_1-1.
 \end{equation}
 Note that $V'$ and $W'$ in \eqref{equation T'=} are allowed to be empty sequences $\varepsilon$.
 By \eqref{equation yij in Mi|Ji|}, we have that
  \begin{equation}\label{equation pi mathop bullet limits_jin 1 t}
 \pi(\mathop{\bullet}\limits_{j\in [1,n]} y_{1, j})\in M_1^{n} \setminus M_1^{n+1}.
 \end{equation}
It follows from \eqref{equation tildeyij=yij}, \eqref{equation tildeyij=1} and \eqref{equation W' is a subsequence of} that for any  $h\in [0,k_1]$,
 \begin{equation}\label{equation the long congurence equation}
\pi(W' \cdot (\mathop{\bullet}\limits_{j\in [1,n]} \widetilde{y}_{1, j}))=\pi(W')*\pi(\mathop{\bullet}\limits_{j\in [1,n]} \widetilde{y}_{1, j})
\equiv 1_R* \pi(\mathop{\bullet}\limits_{j\in [1,n]} y_{1, j})=\pi(\mathop{\bullet}\limits_{j\in [1,n]} y_{1, j}) \pmod {M_1^{h}} .
  \end{equation}
  Since $\pi(V')\in {\rm U}(R)$, it follows from \eqref{equation T'=} that $\pi(T')=\pi(V')*\pi(W' \cdot (\mathop{\bullet}\limits_{j\in [1,n]} \widetilde{y}_{1, j}))$ is an associate of $\pi(W' \cdot (\mathop{\bullet}\limits_{j\in [1,n]} \widetilde{y}_{1, j}))$. Combined with \eqref{equation pi mathop bullet limits_jin 1 t} and \eqref{equation the long congurence equation}, we conclude that
\begin{equation}\label{equation pi(T') in M1tsetminus M1t+1}
\pi(T')\in M_1^{n} \setminus M_1^{n+1}.
\end{equation}
 Since $\pi(T')$ is an idempotent, it follows from \eqref{equation 1 leq t leq k_1-1}  that $\pi(T')=\pi(T')*\pi(T')\in M_1^{2n}\subseteq M_1^{n+1}$, a contradiction with \eqref{equation pi(T') in M1tsetminus M1t+1}. This proves the assertion. \qed

By \eqref{equation |V|_D-1}, \eqref{equation form of the sequence T} and the assertion above, we have ${\rm I}(\mathcal{S}_R)\geq |T|+1={\rm D}({\rm U}(R))+ \sum\limits_{i=1}^r (k_i-1)={\rm D}({\rm U}(R))+ \sum\limits_{i=1}^r ({\rm Ind}(M_i)-1)$ proved. It remains to show the equality ${\rm I}(\mathcal{S}_R)={\rm D}({\rm U}(R))+ \sum\limits_{i=1}^r (k_i-1)$ in the case when $R$ is a local ring or all maximal ideals of $R$ have the indices one, i.e., when $r=1$ or $k_1=\cdots =k_r=1$.

Since $R$ is Artinian, we know that the Jacobson radical $\mathcal{J}(R)=\bigcap\limits_{i=1}^r M_i=\prod\limits_{i=1}^r M_i$ is nilpotent (also as remarked in Section 1 since $R$ is finite then  $M_1,\ldots,M_r$ are all distinct prime ideals of $R$ and $\bigcap\limits_{i=1}^r M_i=\bigcap\limits_{\mbox{\tiny $P$ is a prime ideal of $R$}} P={\rm nil}(R)$ is the nilradical of $R$),
i.e., $(\prod\limits_{i=1}^r M_i)^N={\bf 0}$ for some positive integer $N$. This implies that
\begin{equation}\label{equation product of Miki=0}
\prod\limits_{i=1}^r M_i^{k_i}={\bf 0}
\end{equation}
 is the zero ideal of $R$.

Now we suppose $r=1$, i.e., $M_1$ is the unique maximal ideal.  To prove ${\rm I}(\mathcal{S}_R)={\rm D}({\rm U}(R))+ \sum\limits_{i=1}^r (k_i-1)={\rm D}({\rm U}(R))+ k_1-1,$ take an arbitrary sequence $L$ of terms from $R$ of length ${\rm D}({\rm U}(R))+k_1-1$, it suffices to show that $L$ is not idempotent-product free. Since ${\rm U}(R)=R\setminus M_1$, we have a partition $L=L_1\cdot L_2$ where $L_1$ is a sequence of terms from ${\rm U}(R)$  and $L_2$ is a sequence of terms from $M_1$. By the pigeonhole principle, we have that either (i) $|L_1|\geq {\rm D}({\rm U}(R))$, or (ii)
 $|L_2|\geq k_1$. If (i) holds, then $1_R\in \prod(L_1)\subseteq \prod(L)$, and so $L$ is not idempotent-product free. Otherwise, (ii)
 $|L_2|\geq k_1$ holds, by \eqref{equation product of Miki=0}, then $\pi(L_2)\in M_1^{k_1}={\bf 0}$ which implies that $\pi(L_2)=0_R$ is an idempotent, done.

Suppose $k_1=\cdots =k_r=1$. To prove ${\rm I}(\mathcal{S}_R)={\rm D}({\rm U}(R))+ \sum\limits_{i=1}^r (k_i-1)={\rm D}({\rm U}(R))$, we take an arbitrary sequence $L$ of terms from $R$ and of length ${\rm D}({\rm U}(R))$. It suffices to show that $L$ is not idempotent-product free.
By the Chinese Remainder Theorem,
for any term $a$ of $L$ we can take an element $a'\in R$  such that for each $i\in [1,r]$,
\begin{equation}\label{equaiton tilde a}
a'\equiv\left\{ \begin{array}{ll}
1_R \pmod {M_i}& \textrm{if $a\in M_i$;}\\
a \pmod {M_i} & \textrm{otherwise.}\\
\end{array} \right.
\end{equation}
It follows that $a'\notin M_1\cup \cdots \cup M_r$ and so $a'\in {\rm U}(R)$. Since $|\mathop{\bullet}\limits_{a\mid L} a'|=|L|={\rm D}({\rm U}(R))$, it follows that $1_R\in \prod (\mathop{\bullet}\limits_{a\mid L} a')$, i.e., there exists a nonempty subsequence $W$ of $L$ such that $\prod\limits_{a\mid W}a'=1_R$. Combined with \eqref{equaiton tilde a}, we derive that for each $i\in [1,r]$, either $\pi(W)\equiv 0_R \pmod {M_i}$  or $\pi(W)\equiv \prod\limits_{a\mid W}a'=1_R \pmod {M_i}$, which implies $\pi(W)*\pi(W)\equiv \pi(W)\pmod {M_i}$ in any case. Then $\pi(W)*\pi(W)\equiv \pi(W) \pmod {\bigcap\limits_{i=1}^r M_i}$.  By \eqref{equation product of Miki=0}, $\bigcap\limits_{i=1}^r M_i=\prod\limits_{i=1}^r M_i={\bf 0}$, we have that  $\pi(W)*\pi(W)=\pi(W)$ and so $L$ is not idempotent-product free. This completes the proof of the theorem. \qed

\section{Concluding remarks}

To show that Theorem \ref{Theorem main} implies Theorem E, we shall prove the fact that in Theorem E the quantity $\Omega(K)-\omega(K)$ coincides with
$\sum\limits_{M} ({\rm Ind}(M)-1)$, where $M$ is taken over all distinct maximal ideals of the quotient ring $D\diagup K$. The arguments are as follows.

\noindent {\sl Proof.} \ Since $D$ is a Dedekind domain, then the nonzero proper ideal $K$ has a prime factorization, say $K=P_1^{k_1}*P_2^{k_2}*\cdots *P_r^{k_r}$ where $r>0$, $P_1,P_2,\ldots,P_r$ are distinct prime ideals (maximal ideals) and $k_1,k_2,\ldots, k_r>0$.
Let $\theta: D\rightarrow D\diagup K$ be the canonical epimorphism of $D$ onto the quotient ring $D\diagup K$.
Since $P_1,\ldots, P_r$ are all maximal ideals containing $Ker(\theta)=K$ in the ring $D$ , we know that $\theta(P_1),\theta(P_2),\ldots, \theta(P_r)$ are all distinct maximal ideals of $D\diagup K$.
We see that $P_i\supsetneq P_i^2 \supsetneq \cdots \supsetneq P_i^{k_i}$, $P_i^{t_i}\supseteq K$ and $\theta(P_i)^{t_i}=\theta(P_i^{t_i})$ where $i\in [1,r]$ and $t_i\in [1,k_i]$. It follows that $\theta(P_i)\supsetneq \theta(P_i)^2 \supsetneq \cdots \supsetneq \theta(P_i)^{k_i}$, and so
\begin{equation}\label{equation ind Pigeq ki}
{\rm Ind}(\theta(P_i))\geq k_i \mbox{ for each } i\in [1,r].
\end{equation}
On the other hand, for each $i\in [1,r]$ we have that $\theta^{-1}(\theta(P_i^{k_i+1}))=P_i^{k_i+1}+K=P_i^{k_i+1}+\prod\limits_{j\in [1,r]} P_j^{k_j}=P_i^{k_i}*(P_i+\prod\limits_{j\in [1,r]\setminus \{i\}} P_j^{k_j})=P_i^{k_i} * R=P_i^{k_i}=P_i^{k_i}+\prod\limits_{j\in [1,r]} P_j^{k_j} = P_i^{k_i}+K = \theta^{-1}(\theta(P_i^{k_i}))$. Since $\theta$ is surjective, we derive that $\theta(P_i^{k_i+1})=\theta(P_i^{k_i})$ and so $\theta(P_i)^{k_i+1}=\theta(P_i)^{k_i}$. Combined with \eqref{equation ind Pigeq ki}, we have that ${\rm Ind}(\theta(P_i))=k_i \mbox{ for each } i\in [1,r].$
Then $\Omega(K)-\omega(K)=\sum\limits_{i=1}^r k_i-r=\sum\limits_{i=1}^r ({\rm Ind}(\theta(P_i))-1)=\sum\limits_{M} ({\rm Ind}(M)-1)$ where $M$ is taken over all distinct maximal ideals of $D\diagup K$,  done. \qed

We close this paper with the following problem.

\noindent {\bf Problem.}  Let $R$ be a finite commutative unitary ring.  Determine when ${\rm I}(\mathcal{S}_R)={\rm D}({\rm U}(R))+\sum\limits_{M} ({\rm Ind}(M)-1)$ where $M$ is taken over all distinct maximal ideals of $R$.

\end{document}